\documentclass[11pt,a4paper]{article}
\usepackage{amsfonts}
\newtheorem{theo}{Theorem}
\newenvironment{Theo}{\begin{theo}\slshape}{\end{theo}}
\newtheorem{Lemm}{Lemma}[section]

\newtheorem{Corol}{Corollary}[section]
\newtheorem{rema}{Remark}

\newcommand{\p}{\partial}

\def\qed{\hfill$\square$\par \bigskip}
\newenvironment{Demo}[1]{{\bf Proof #1.~}}{\qed}
\newcommand{\R}{\mathbb{R}}

\newcommand{\para}[1]{\left(#1\right)}

\usepackage{graphicx}

\usepackage{mathrsfs}
\usepackage[latin1]{inputenc}
\usepackage[T1]{fontenc}
\usepackage{color,xspace}
\usepackage{hyperref}
\usepackage[usenames,dvipsnames]{pstricks}
\usepackage{times}

\begin{document}
\title{Stability estimate for hyperbolic inverse problem with time dependent coefficient}
\author{
\small {\bf Ibtissem Ben A\"{\i}cha }\\
\small Department of Mathematics,\\
\small Faculty of Sciences of Bizerte,\\
\small 7021 Jarzouna Bizerte, Tunisia. \footnote{e-mail correspondence: mourad.bellassoued@fsb.rnu.tn}}
\date{}
\maketitle
\begin{abstract}
 We study the  stability  in the inverse problem of determining the time
dependent zeroth-order coefficient $q(t,x)$ arising in the wave equation,
 from boundary
observations. We derive, in dimension $n\!\geq \!\!2$ , a log-type
stability estimate in  the determination of $q$ from the Dirichlet-to-Neumann map, in a subset of our domain assuming that it is known outside this subset. Moreover, we prove that we can extend this result to the determination of $q$ in a larger region, and then in the whole domain provided that we have much more data. \\
\textbf{Keywords:} Inverse problems,
Dirichlet-to-Neumann map, Wave equation, Bounded domain, Time dependent
potential,  X-ray transform, Stability estimate.
\end{abstract}
\section{Introduction}\label{Sec1}
\subsection{Statement of the problem}\label{Sub1.1}
Let $\Omega$ be a bounded domain of $\R^{n}$, $n\geq 2$, with
$\mathcal{C}^{\infty}$ boundary $\Gamma\!=\!\p\Omega.$ Given $T\!>\!2\,\mbox{Diam}(\Omega)$, we
introduce the following initial boundary value problem for
the wave equation
 \begin{equation}\label{Eq1.1}
\left\{
   \begin{array}{ll}
     (\p_{t}^{2}-\Delta +q(t,x))u=0 & \mbox{in}\,\,Q=[0,T]\times\Omega, \\

 u(0,x)=u_{0},\,\,\p_{t}u(0,x)=u_{1} & \mbox{in}\,\,\Omega, \\

   u=f & \mbox{on}\,\,\Sigma=[0,T]\times
\Gamma,
   \end{array}
 \right.
\end{equation}
where $f\in H^{1}(\Sigma)$, $u_{0}\in H^{1}(\Omega),$ $u_{1}\in L^{2}(\Omega)$ and the potential $q \in\mathcal{C}^{1}(\overline{Q})$ is assumed to be
real valued. It is well-known (see \cite{[C23]}, \cite{[R]}) that if the compatibility condition is satisfied,
then (\ref{Eq1.1}) is well-posed. Therefore we can introduce the following operator
$$\begin{array}{ccc}
\Lambda_{q}:H^{1}(\Sigma)&\longrightarrow& L^{2}(\Sigma)\\
f&\longmapsto&\p_{\nu}u,
\end{array}
$$
usually called the Dirichlet-to-Neumann map. Here $\nu(x)$ denotes the unit
outward normal to $\Gamma$ at $x$ and $\p_{\nu}u$ stands for $\nabla u.\nu.$\\

  In the present paper, we will first study the inverse problem of recovering the time dependent potential $q$ from the Dirichlet-to-Neumann map $\Lambda_{q}$ associated to the problem (\ref{Eq1.1}) with $(u_{0},u_{1})=(0,0)$.
 This inverse problem is to know whether the knowledge of $\Lambda_{q}$,
 can uniquely determine the electric time dependent  potential $q$. \\

 Physically, it consists in determining physical properties
such as the time evolving density of an inhomogeneous medium by probing it
with disturbances generated on the boundary. 
 And the goal is to
recover $q$ which describes the property of the medium. We assume that
the medium is quiet initially and the Dirichlet data $f$ is a disturbance
used to probe it.\\

The problem of recovering coefficients for hyperbolic equations from boundary
measurements was treated by many authors. In \cite{[B3]} Rakesh and Symes
proved  a uniqueness result in recovering the time independent potential
$q(x)$ in the wave equation. In \cite{[B6]} Isakov treated the inverse
problem of recovering a zeroth order coefficient and a damping coefficient.
These results are concerned in the case where  the Dirichlet-to-Neumann map
is considered in the whole boundary. A key ingredient in the existing
results, is the construction of complex geometric optics solutions
concentrating near lines with any direction $\omega\in \mathbb{S}^{n-1}$ and the
relationship  between the hyperbolic Dirichlet-to-Neumann map and the X-ray
transform plays a crucial role. The uniqueness in the determination of time
independent potential appearing in the wave equation by a local observations   was
proved by Eskin \cite{[B7]}.\\

The uniqueness by local measurements is solved well. However, the stability by a local Dirichlet-to-Neumann map is not discussed comprehensively. For it, one can see Bellassoued, Chouli and Yamamoto \cite{[B9]} where a log-type stability estimate was proved in the case where the Neumann data are observed in an arbitrary subdomain of the boundary,   Isakov
and Sun \cite{[B8]} where a local Dirichlet-to-Neumann map yields an H\"{o}lder stability result in determining a coefficient in a subdomain. The case where the Neumann data are observed in the whole boundary, a stability of H\"{o}lder type was established in Cipolatti and Lopez \cite{[B10]}, Sun \cite{[B12]}, and in Riemannian case in
 M. Bellassoued and D. Dos  Santos Ferreira \cite{[B13]}, Stefanov and Uhlmann \cite{[B11]}.\\

 All the above mentioned results are concerned only with time-independent coefficients. Many authors considered the problem of determining time-dependent coefficients for hyperbolic equations. In \cite{[B14]}, Stefanov proved that the time dependent potential $q(t, x)$ arising in the wave equation  is uniquely determined
 from the knowledge of scattering data. In \cite{[B15]}, Ramm and Sj\"{o}strand  treated the problem of determining the time-dependent potential $q(t,x)$ from Dirichlet-to-Neumann map, on the infinite time-space cylindrical domain $\R_{t}\times\Omega$, and they proved a uniqueness result under suitable  assumptions. In \cite{[B16]}, R. Salazar, extended  the results in \cite{[B15]} to more general coefficients and  proved a result of stability for compactly supported coefficients provided $T$ is sufficiently large. \\

The inverse problem of determining the time-dependent coefficient $q(t,x)$ from the Dirichlet-to-Neumann map $\Lambda_{q}$, was treated by Ramm and Rakesh \cite{[B5]}, they assumed without loss of generality that $\Omega$ is a ball and they proved  a uniqueness result only in a subset  made of lines making $45^{\circ}$
with the $t$-axis and meeting the planes $t=0$ and $t=T$ outside
$\overline{Q}$, provided that it's known outside this subset.  It's clear that with zero initial data one can not hope to recover $q(t,x)$ over the whole domain $Q$, even from the knowledge of the full boundary operator $\Lambda_{q}$. This is due to the domain of dependence  associated to the hyperbolic problem (\ref{Eq1.1}) (see \cite{[C2]}).
However, in Isakov \cite{[B17]}, the ideas from \cite{[B18]}-\cite{[B19]} are used to prove a  uniqueness result in determining $q(t,x)$ over the whole domain $Q$, but he needed much more information. Indeed his data was the response of the medium for all possible initial data.\\

In this paper, we will prove a log-type stability estimate  which establishes that
the time dependent potential $q(t,x)$ depends stably on the Dirichlet-to Neumann map $\Lambda_{q}$ in  a subset of our domain, provided that it is known outside this subset. After that we prove that we can  extend this result to the determination of $q$ in a larger region if we further know the measures $(u(T,.),\p_{t}u(T,.))$, where $u$ is the solution of the initial boundary value problem (\ref{Eq1.1}) with $(u_{0},u_{1})=(0,0)$. Moreover, we will prove that if our data was the response of the medium for all possible initial data, then we have a log-type stability estimate for this problem  over the whole domain $Q$.\\

  Inspired by the work of M. Bellassoued and D. Dos Santos Ferreira \cite{[B13]}, Alden Waters \cite{[C3]} succeeded in proving  a type of an H\"{o}lder stability estimate for the inverse problem of recovering the X-ray transform of the time-dependent potential $q$, appearing in the wave equation, from the dynamical Dicrichlet-to Neumann map in Riemannian case. A key ingredient in this result is the construction of Gaussian beam solutions. In the case $n\geq 3$, the inverse problem associated to the system (\ref{Eq1.1}) with the initial condition $u_{0}=0$, was treated recently by Y. Kian \cite{[C24]}, indeed, inspired by Bellassoued-Jellali-Yamammoto \cite{[B2]}-\cite{[Corr]} and using suitable complex geometric optics solutions and Carleman estimate, he proved a log-log type stability estimate in determining the time dependent coefficient $q(t,x)$, from the knowledge of  partial Dirichlet-to-Neumann measurement and the measure $u(T,.)$.\\

Before stating our main results, we recall the following Lemma on the unique existence of a solution to the problem (\ref{Eq1.1}). The proof is given in \cite{[C23]} (see also \cite{[R]}).
\begin{Lemm}\label{Lm1.1}
Let $T>0$ be given. Suppose that
$u_{0}\in H^{1}(\Omega),\,\,\,u_{1}\in L^{2}(\Omega),\,\,\,\,\mbox{and}\,\,f\in H^{1}(\Sigma).$
Assume, in addition, that  $f(0,.)=u_{0}|_{\Gamma}.$
Then, there exists a unique solution $u$ of (\ref{Eq1.1}) satisfying
 $$u\in \mathcal{C}([0,T];\,H^{1}(\Omega))\cap\mathcal{C}^{1}([0,T];\,L^{2}(\Omega)),$$
and there exists $C>0$ such that for any $t\in [0,T]$, we have
$$\begin{array}{lll}
\|\p_{\nu}u\|_{L^{2}(\Sigma)}+\|u(t,.)\|_{H^{1}(\Omega)}+\|\p_{t}u(t,.)\|_{L^{2}(\Omega)}
\!&\leq&\!\!C\displaystyle\para{\|f\|_{H^{1}(\Sigma)}+\|u_{0}\|_{H^{1}(\Omega)}+\|u_{1}\|_{L^{2}(\Omega)}}.
\end{array}$$
\end{Lemm}
From the above Lemma one can see that, if $(u_{0},u_{1})=(0,0)$, the Dirichlet-to-Neumann map $\Lambda_{q}$ is continuous from $H^{1}(\Sigma)$ to $L^{2}(\Sigma)$. Therefore we denote by $\|\Lambda_{q}\|$ its norm in $\mathcal{L}(H^{1}(\Sigma),\,L^{2}(\Sigma))$.\\
\subsection{Main results}\label{Sub1.2}
In order to state our main results  we first introduce some notations:\\
\\
Let $r>0$ such that $T>2r$ and $\overline{\Omega}\subseteq B(0,\frac{r}{2})=\displaystyle\left\{x\in\R^{n},\,|x|\leq\displaystyle\frac{r}{2}\right\}$. We set  $Q_{r}=[0,T]\times B(0,\frac{r}{2}).$
We consider the following sets

$$\mathscr{A}_{r}=\left\{x\in\R ^{n},\,\,\displaystyle\frac{r}{2}<|x|<T-\displaystyle\frac{r}{2}\right\}.$$
$$
\mathscr{C}_{r}^{+}=\left\{(t,x)\in Q_{r},\,\,|x|<t-\displaystyle\frac{r}{2},\,t>\frac{r}{2}\right\}.$$
$$\mathscr{C}_{r}^{-}=\left\{(t,x)\in Q_{r},
\,\,|x|<T-\displaystyle\frac{r}{2}-t,\,T-\frac{r}{2}> t\right\}.$$

 Note also $Q_{r}^{*}=\mathscr{C}_{r}^{+}\cap \mathscr{C}_{r}^{-}$. Let denote by  $Q_{*}= Q\cap Q_{r}^{*}$. We remark that $Q_{*}$ is made of lines making $45^{\circ}$
 with the $t$-axis and meeting the planes $t=0$ and $t=T$ outside
$\overline{Q}_{r}$.  We denote by  $Q_{\sharp}=Q\cap \mathscr{C}_{r}^{+}$. We remark that $Q_{\sharp}$ is made of lines making $45^{\circ}$
with the $t$-axis and meeting only the planes $t=0$ outside $\overline{Q}_{r}$. Let's note that $Q_{*}\subset Q_{\sharp}\subset Q$.
\begin{rema}
In the particular case where $\overline{\Omega}=B(0,\frac{r}{2}),$ we remark that $Q_{*}=Q_{r}^{*}$ which is  the region $I$ in Figure\ref{Fig}. And $Q_{\sharp}=\mathscr{C}_{r}^{+}$ which is the region $I\,\cup\, II\,\cup\,III\,\cup\,IV$.
\end{rema}
\begin{center}\label{Fig}
\includegraphics[width=16cm,height=11cm]{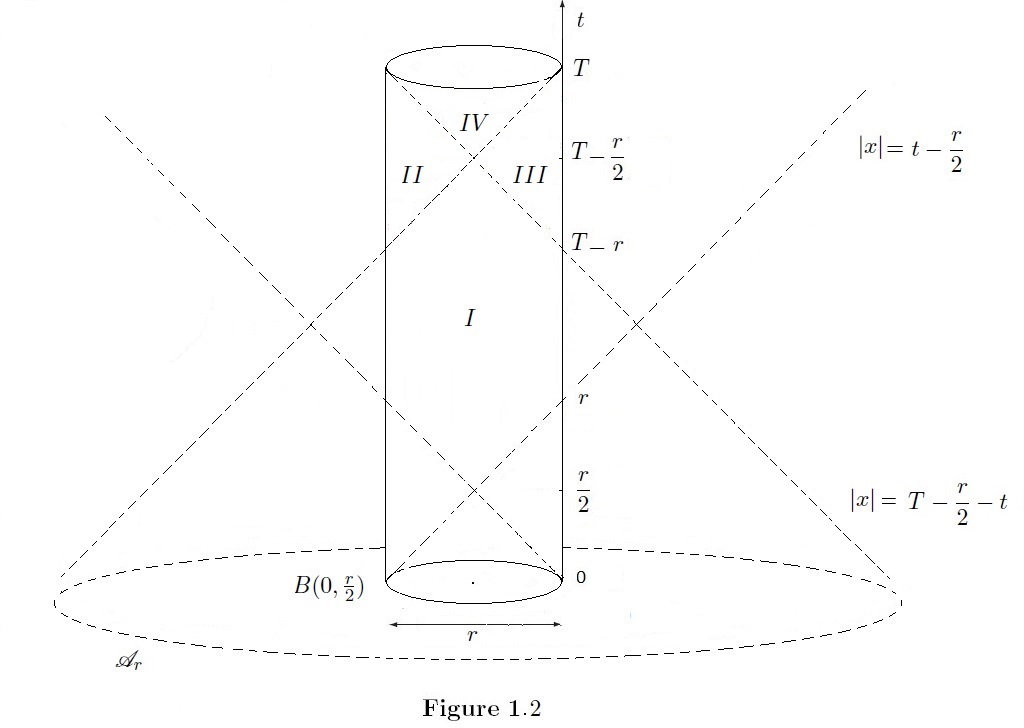}
\end{center}
Further, given $q_{0}\in\mathcal{C}^{1}(\overline{Q}_{r})$ and $M>0$, we introduce
$$\mathcal{A}^{*}(q_{0},M)=\left\{q\in \mathcal{C}^{1}(\overline{Q}_{r}),\,\,q=q_{0}\,\,\mbox{in}\,\,\overline{Q}_{r}\setminus Q_{*},\,\,\|q\|_{L^{\infty}(Q)}\leq M\right\},$$
and
$$\mathcal{A}^{\sharp}(q_{0},M)=\left\{q\in \mathcal{C}^{1}(\overline{Q}_{r}),\,\,q=q_{0}\,\,\mbox{in}\,\,\overline{Q}_{r}\setminus  Q_{\sharp},\,\,\|q\|_{L^{\infty}(Q)}\leq M\right\}.$$
Then our first main result can be stated as  follows:
\begin{Theo}\label{Thm1}
Assume that $T>2\,\mbox{Diam}\,(\Omega)$. Then, for every $q_{1},\,q_{2}\in \mathcal{A}^{*}(q_{0},M)$, there
exist two constants $C>0$ and $\mu_{1}\in (0,1)$, such that we have
$$\|q_{1}-q_{2}\|_{H^{-1}(Q_{*})}\leq
C\,\para{\|\Lambda_{q_{1}}-\Lambda_{q_{2}}\|^{\mu_{1}}+|\log\|\Lambda_{q_{1}}-\Lambda_{q_{2}}\|
|^{-1}},$$ where $C$ depends only on $\Omega,\,M,\,T,$ and $n$.

Suppose  in addition that $q_{1},\,q_{2}\in H^{s+1}(Q),$ for
$s>\displaystyle\frac{n}{2}$ and that $\|q_{i}\|_{H^{s+1}(Q)}\leq
M,$ $i=1,2$, for some $M>0$, then there exist two constants $C'>0$ and
$\mu_{2}\in(0,1)$ such that
\begin{equation}\label{Eq1.2}
\|q_{1}-q_{2}\|_{L^{\infty}(Q_{*})}\leq C'\para{\|\Lambda_{q_{1}}-\Lambda_{q_{2}}\|+|\log\|\Lambda_{q_{1}}-\Lambda_{q_{2}}\||^{-1}}^{\mu_{2}}.\end{equation}

\end{Theo}
As an immediate consequence of Theorem \ref{Thm1}, we have the following uniqueness
result.
\begin{Corol} (Uniqueness)
Under the same assumptions, for every $q_{1},\,q_{2}\in\mathcal{A}^{*}(q_{0},M),$ we have
the uniqueness
$$\Lambda_{q_{1}}(f)=\Lambda_{q_{1}}(f),\,\,\mbox{for}\,\,\mbox{any}\,f\in H^{1}(\Sigma),\,\,\,\mbox{imply}\,\,\,\,q_{1}(t,x)=q_{2}(t,x),$$
everywhere in $Q_{*}$.
\end{Corol}

Let us note that in this result we determine the time dependent coefficient $q$ from full boundary measurements $\Lambda_{q}$ only in a subset $Q_{*} \subset Q$, provided that it is known outside of this part. \\

In order to extend this result to the determination of $q$ in a larger region $Q_{\sharp}\supset Q_{*}$ we need more information about the solution $u$. Namely we need the measures of $(u(T,.),\p_{t}u(T,.))$.
So, let's introduce the following  boundary operator:
$$\begin{array}{ccc}
\mathscr{R}_{q}:H^{1}(\Sigma)&\longrightarrow& L^{2}(\Sigma)\times H^{1}(\Omega)\times L^{2}(\Omega).\\
f&\longmapsto&(\p_{\nu}u,\,u(T,.),\,\p_{t}u(T,.))
\end{array}
$$
From Lemma \ref{Lm1.1}, we deduce that, if $(u_{0},u_{1})=(0,0)$, the operator $\mathscr{R}_{q}$ is continuous from $H^{1}(\Sigma)$ to $L^{2}(\Sigma)\times H^{1}(\Omega)\times L^{2}(\Omega)$. We denote by $\|\mathscr{R}_{q}\|$ its norm in $\mathcal{L}\Big(H^{1}(\Sigma),\,L^{2}(\Sigma)\times H^{1}(\Omega)\times L^{2}(\Omega)\Big).$
\\
\\
Then, the second result is the following:
\begin{Theo}\label{Thm2}
Assume that $T>2\,\mbox{Diam}\,(\Omega)$. Then, for every $q_{1},\,q_{2}\in \mathcal{A}^{\sharp}(q_{0},M)$, there
exist two constants $C>0$ and $\mu_{1}\in (0,1)$, such that we have
$$\|q_{1}-q_{2}\|_{H^{-1}(Q_{\sharp})}\leq
C\,\para{\|\mathscr{R}_{q_{1}}-\mathscr{R}_{q_{2}}\|^{\mu_{1}}+|\,\, \log\|\mathscr{R}_{q_{1}}-\mathscr{R}_{q_{2}}\|
|^{-1}},$$ where $C$ depends only on $\Omega,\,M,\,T,$ and $n$.

Suppose  in addition that $q_{1},\,q_{2}\in H^{s+1}(Q),$ for
$s>\displaystyle\frac{n}{2}$ and that $\|q_{i}\|_{H^{s+1}(Q)}\leq
M,$ $i=1,2$, for some $M>0$, then there exist two constants $C'>0$ and
$\mu_{2}\in(0,1)$ such that
$$\|q_{1}-q_{2}\|_{L^{\infty}(Q_{\sharp})}\leq C'\para{\|\mathscr{R}_{q_{1}}-\mathscr{R}_{q_{2}}\|+|\ \log\|\mathscr{R}_{q_{1}}-\mathscr{R}_{q_{2}}\||^{-1}}^{\mu_{2}}.$$
where $C'$ depends on $\Omega$, $M$, $T$, and $n$.
\end{Theo}
As an immediate consequence of Theorem \ref{Thm2}, we have the following uniqueness
result.
\begin{Corol} (Uniqueness)
Under the same assumptions, for every $q_{1},\,q_{2}\in\mathcal{A}^{\sharp}(q_{0},M),$ we have
the uniqueness
$$\mathscr{R}_{q_{1}}(f)=\mathscr{R}_{q_{2}}(f),\,\,\mbox{for\,\,any}\,\,f\in H^{1}(\Sigma),\,\,\,\mbox{imply}\,\,\,\,q_{1}(t,x)=q_{2}(t,x),$$
everywhere in $Q_{\sharp}$.
\end{Corol}

With zero initial data there is no hope to recover $q(t,x)$ over the whole domain $Q$,  even from  the knowledge of the boundary operator $\mathscr{R}_{q}$. However, from measurements made for all possible initial data, we can extend the results in Theorem \ref{Thm1} and Theorem \ref{Thm2} to the determination of $q$ over the whole domain. We define  the boundary operator
$$\begin{array}{ccc}
\mathcal{I}_{q}: H^{1}(\Sigma)\times H^{1}(\Omega)\times L^{2}(\Omega)&\longrightarrow& L^{2}(\Sigma)\times H^{1}(\Omega)\times L^{2}(\Omega).\\
(f,u_{0},u_{1})&\longmapsto&(\p_{\nu}u,\,u(T,.),\,\p_{t}u(T,.))
\end{array}$$
From Lemme \ref{Lm1.1}, we deduce that the linear operator $\mathcal{I}_{q}$ is continuous from $H^{1}(\Sigma)\!\times \!H^{1}(\Omega)\!\times\! L^{2}(\Omega)$ to $L^{2}(\Sigma)\!\times\! H^{1}(\Omega)\!\times\! L^{2}(\Omega)$. We denote by $\|\mathcal{I}_{q}\|$ its norm.\\
\\
Then, our last result can be stated as follows:

\begin{Theo}\label{Thm3}
Assume that $T>2\,\mbox{Diam}\,(\Omega)$. Then, for every $q_{1},\,q_{2}\in \mathcal{C}^{1}(\overline{Q})$, such that $\|q_{i}\|_{L^{\infty}(Q)}\leq M,$ for $i=1,\,2$. There exist two constants $C>0$ and $\mu_{1}\in(0,1)$, such that we have $$
\|q_{1}-q_{2}\|_{H^{-1}(Q)}\leq C\para{ \|\mathcal{I} _{q_{1}}-\mathcal{I}_{q_{2}}\|^{\mu_{1}}+| \log \|\mathcal{I}_{q_{1}}-\mathcal{I}_{q_{2}}\||^{-1}},
$$
where $C$ depends only on $\Omega, \,\,M\,,\,T,\,\,\mbox{and}\,\,n$.

Suppose in addition that $q_{1},\,q_{2}\in H^{s+1}(Q)$, for $s>\frac{n}{2}$ and $\|q_{i}\|_{H^{s+1}(Q)}\leq M,$ $i=1,\,2$, for some $M>0$, then
there exist two constants $C'>0$ and $\mu_{2}\in(0,1)$ such that
$$
\|q_{1}-q_{2}\|_{L^{\infty}(Q)}\leq  C' \para{\|\mathcal{I}_{q_{1}}-\mathcal{I}_{q_{2}}\|+| \log\|\mathcal{I}_{q_{1}}-\mathcal{I}_{q_{2}}\||^{-1}}^{\mu_{2}}.
$$
 \end{Theo}
 As an immediate consequence of Theorem \ref{Thm3}, we have:
\begin{Corol}
Under the same assumptions as in Theorem \ref{Thm3}, we have the uniqueness
$$\mathcal{I}_{q_{1}}=\mathcal{I}_{q_{2}},\,\,\mbox{imply}\,\, q_{1}(t,x)=q_{2}(t,x), \,\,\,\mbox{in}\,\,Q.$$
\end{Corol}

This paper is organized as follows. In section \ref{Sec2} we construct special optics geometrical solutions to the wave equation (\ref{Eq1.1}).
Using these geometric optics solutions, in section \ref{Sec3} we prove Theorem \ref{Thm1}, in section \ref{Sec4} we prove Theorem \ref{Thm2} and in section \ref{Sec5} we prove Theorem \ref{Thm3}.

\section{Geometric optics solutions}\label{Sec2}
In the present section, we collect some results which are needed in the proof of our main results. We start by the following Lemma (see \cite{[C23]}, \cite{[R]}):
\begin{Lemm}\label{Lm2.1}
Let $T>0$ and $q\in L^{\infty}(Q)$, suppose that $F\in L^{1}(0,T;L^{2}(\Omega))$. The unique solution $u$ of the system
$$
\left\{
  \begin{array}{ll}
   \para{\p_{t}^{2}-\Delta+q(t,x)} u(t,x)=F(t,x)  & \mbox{in}\,\, Q, \\
    u(0,x)=\p_{t}u(0,x)=0 & \mbox{in}\,\,\Omega, \\
    u(t,x)=0 & \mbox{on}\,\,\Sigma,
  \end{array}
\right.
$$
satisfies $$u\in \mathcal{C}([0,T];H^{1}_{0}(\Omega))\cap\mathcal{C}^{1}([0,T];L^{2}(\Omega)).$$
Moreover, there exists a constant $C>0$ such that
\begin{equation}\label{Energy}
\|\p_{t}u(t,.)\|_{L^{2}(\Omega)}+\|\nabla u(t,.)\|_{L^{2}(\Omega)}\leq C\|F\|_{L^{1}(0,T;L^{2}(\Omega))}.
\end{equation}
\end{Lemm}
Using Lemma \ref{Lm2.1}  we are able to construct suitable geometrical optics solutions for
our inverse problem, which are key ingredients to the  proof of  our main
results.\\
\\
Let $\varphi\in \mathcal{C}_{0}^{\infty}(\R^{n})$. Notice that for all $\omega\in\mathbb{S}^{n-1}=\{\omega\in\R^{n},\,\,|\omega|=1\},$ the function
\begin{equation}\label{Eq2.3}
a(t,x)=\varphi(x+t\omega)\end{equation}
solves the transport equation
\begin{equation}\label{Eq2.4}
(\p_{t}-\omega.\nabla)a(t,x)=0.
\end{equation}
Let's now prove the following Lemma:
\begin{Lemm}\label{Lem2.2}
Let $q\in\mathcal{C}^{1}(\overline{Q})$ such that $\|q\|_{L^{\infty}(Q)}\leq M$.
For $\omega\in \mathbb{S}^{n-1},$ and  $\varphi\in\mathcal{C}_{0}^{\infty}(\R^{n})$, we consider the function $a$ defined
by (\ref{Eq2.3}). Then, for $\lambda>0$, the equation
\begin{equation}\label{Eq2.5}
\para{\p_{t}^{2}-\Delta+q(t,x)}u(t,x)=0\,\,\mbox{in} \,\,Q,
\end{equation}
admits a solution
$$u^{\pm}\in\mathcal{C}([0,T];H^{1}(\Omega))\cap\mathcal{C}^{1}([0,T];L^{2}(\Omega)),$$
of the following form
\begin{equation}\label{Eq2.6}
u^{\pm}(t,x)=a(t,x) e^{\pm \,i\,\lambda(x.\omega+t)}+R^{\pm}(t,x),\end{equation}
where $R^{\pm}(t,x)$ satisfies
$$R^{\pm}(t,x)=0,\,\,\,\mbox{for\,all}\,\,(t,x)\in\Sigma$$
and
$$\p_{t}R^{+}(0,x)=R^{+}(0,x)=0,\,\,\,x\in\Omega,$$
$$\p_{t}R^{-}(T,x)=R^{-}(T,x)=0,\,\,\,x\in\Omega.$$
Moreover,
\begin{equation}\label{Eq2.7}
\lambda\,\|R^{\pm}\|_{L^{2}(Q)}+\|\nabla R^{\pm}\|_{L^{2}(Q)}\leq C\,\|\varphi\|_{H^{3}(\R^{n})},
\end{equation}
where $C$ depends only on $\Omega,\,\,\,T\,\,\mbox{and}\,\,\,M.$
\end{Lemm}
\begin{Demo}{}
We adapt the strategy developed in the proof of  a similar result in
\cite{[B3]} , where a time independent potential $q$ was considered. In light
of (\ref{Eq2.5}) and (\ref{Eq2.6}) it is enough to prove the existence of
$R^{\pm}$ satisfying
\begin{equation}\label{Eq2.8}
\left\{
  \begin{array}{ll}
    \displaystyle\para{\p^{2}_{t}-\!\!\Delta+q(t,x)}R^{\pm}(t,x)\!\!=\!\!-\displaystyle\para{\p_{t}^{2}-\Delta+q(t,x)}\displaystyle\para{a(t,x) e^{\pm i \,\lambda\,(x.\omega+t)}}  & \mbox{in}\,\, Q, \\
    R^{\pm}(\theta,x)=0,\,\,\p_{t}R^{\pm}(\theta,x)=0 ,\,\,\,\,\,\,\,\,\,\,\theta=0,\,\,\mbox{or}\,\,\,T& \mbox{in}\,\,\Omega, \\
  R^{\pm}(t,x)=0 & \mbox{on}\,\,\Sigma,
  \end{array}
\right.
\end{equation}
and obeying (\ref{Eq2.7}). We prove the result for  $u^{+}$. The existence of $u^{-}$, being handled in a similar way.
To do that note
$$
g(t,x)=-\para{\p_{t}^{2}-\Delta+q(t,x)}\para{a(t,x) e^{i \,\lambda\,(x.\omega+t)}}
$$
and use  (\ref{Eq2.4}), getting
\begin{equation}\label{Eq2.9}
g(t,x)
=-e^{i\,\lambda\,(x.\omega+t)}\para{\p_{t}^{2}-\Delta+q(t,x)}a(t,x)=-e^{i\,\lambda\,(x.\omega+t)}g_{0}(t,x),
\end{equation}
where $g_{0}\in L^{1}(0,T; L^{2}(\Omega))$. Thus, $R$  is a suitable solution to the system (\ref{Eq2.8}) satisfying
$$R\in\mathcal{C}([0,T];H^{1}_{0}(\Omega))\cap\mathcal{C}^{1}([0,T];L^{2}(\Omega))$$
and the function
\begin{equation}\label{Eq2.10}
w(t,x)=\int_{0}^{t}R(s,x)\,ds
\end{equation}
solves the following equation
$$
\left\{
  \begin{array}{ll}
    \para{\p_{t}^{2}-\Delta+q(t,x)} w(t,x)=F_{1}(t,x)+F_{2}(t,x) & \mbox{in}\,\,Q, \\
    w(0,x)=0,\,\,\p_{t}w(0,x)=0 & \mbox{in}\,\,\Omega, \\
    w(t,x)=0 & \mbox{on}\,\,\Sigma.
  \end{array}
\right.
$$
Where
\begin{equation}\label{Eq2.11}
F_{1}(t,x)=\displaystyle\int_{0}^{t}g(s,x)\,ds,\,\,\,\mbox{and}\,\,\,F_{2}(t,x)=\displaystyle\int_{0}^{t}[q(t,x)-q(s,x)]R(s,x)\,ds.
\end{equation}
Let $\tau\in [0,T]$. In use of Lemma \ref{Lm2.1} on the interval $[0,\tau]$, there exists a constant $C>0$ such that
\begin{eqnarray}\label{Eq2.12}
\|\p_{t}w(\tau,.)\|^{2}_{L^{2}(\Omega)}&\leq& C \|F_{1}+F_{2}\|^{2}_{L^{1}(0,\tau;L^{2}(\Omega))}\cr
&\leq& C \displaystyle\para{\|F_{1}\|^{2}_{L^{2}(Q)}+\|F_{2}\|^{2}_{L^{2}(0,\tau;L^{2}(\Omega))}}.
\end{eqnarray}
Using (\ref{Eq2.10}), we have
$$\begin{array}{lll}
\|F_{2}\|^{2}_{L^{2}(0,\tau,L^{2}(\Omega))}&\leq& C_{T}\|q\|_{L^{\infty}(Q)}^{2}\,\displaystyle\int_{0}^{\tau}\|\p_{t}w(s,.)\|_{L^{2}(\Omega)}^{2}ds.
\end{array}$$
Then, it follows from (\ref{Eq2.12}) that
$$
\|\p_{t}w(\tau,.)\|^{2}_{L^{2}(\Omega)}\leq C\para{\|F_{1}\|^{2}_{L^{2}(Q)}+\|q\|^{2}_{L^{\infty}(Q)}\int_{0}^{\tau}\|\p_{t}w(s,.)\|^{2}_{L^{2}(\Omega)}\,ds}.
$$
Then, from Gronwall's inequality, one gets
$$
\|\p_{t}w(\tau,.)\|^{2}_{L^{2}(\Omega)}\leq C_{T}\,\|F_{1}\|^{2}_{L^{2}(Q)},\,\,\,\,
$$
 where the constant $C_{T}>0$ depends on $T$ and $\|q\|_{L^{\infty}}.$  From where we get
\begin{eqnarray}\label{Eq2.13}
\|R\|^ {2}_{L^{2}(Q)}\leq C_{T}\,\|F_{1}\|^{2}_{L^{2}(Q)},\,
\end{eqnarray}
\mbox{according\, to \,}(\ref{Eq2.10}). Further, as
$$
\|F_{1}\|_{L^{2}(Q)}^{2}=\frac{1}{\lambda^{2}}\int_{Q}|\int_{0}^{t}g_{0}(s,x)\,\p_{s}(e^{i\lambda(x.\omega+s)})\,ds|^{2}\,dx\,dt,\,\,\,
$$
\mbox{by\,}(\ref{Eq2.9}) \,\mbox{and}\,(\ref{Eq2.11}). Then, integrating by parts with respect to $s$, we deduce from $(\ref{Eq2.13})$ that
there exists a constant $C>0$ such that
$$
\|R\|_{L^{2}(Q)}\leq \frac{C}{\lambda}\,\|\varphi\|_{H^{3}(\R^{n})}.
$$
Finally, Since $\|g\|_{L^{2}(Q)}\leq C\,\|\varphi\|_{H^{3}(\R^{n})}$, using the energy estimate (\ref{Energy}) for the problem (\ref{Eq2.8}) we obtain
$$\|\nabla R\|_{L^{2}(Q)}\leq C_{T}\,\|\varphi\|_{H^{3}(\R^{n})},$$
 This completes the proof.
\end{Demo}
\section{Proof of Theorem \ref{Thm1}}\label{Sec3}
In the present section we will prove a log-type stability estimate in determining $q$ appearing in the initial boundary value problem (\ref{Eq1.1}) with $(u_{0},u_{1})=(0,0)$.
The main ingredients of the proof are
geometric optics solutions introduced in Section \ref{Sec2} and X-ray transform. We start by considering geometric optics solutions of the form (\ref{Eq2.6}). We only assume that supp$\,\varphi\subset \mathscr{A}_{r}$,  in such a way we have
$$\mbox{supp}\,\varphi\cap\Omega=\emptyset,\,\,\,\mbox{and}\,\,\,\para{\mbox{supp}\,\varphi\pm T\omega}\cap\Omega=\emptyset,\,\,\forall\,\omega\in \mathbb{S}^{n-1}.$$
Then we have the following preliminary  estimate which relates the differential of  two potentials to the Dirichlet-to-Neumann map.
\begin{Lemm}\label{Lem3.1}
 Let $q_{1},q_{2}\in \mathcal{A}^{*}(q_{0},M)$, and put $q=\para{q_{2}-q_{1}}$. There exists $C>0,$ such that for any $\omega\in \mathbb{S}^{n-1}$ and
$\varphi\in \mathcal{C}^{\infty}_{0}(\mathscr{A}_{r})$, the following estimate
\begin{equation}\label{Eq3.14}
|\displaystyle\int_{0}^{T}\int_{\R^{n}}q(t,x-t\omega)\,\varphi^{2}(x)\,dx\,dt|\leq C\,\para{\lambda^{3}\|\Lambda_{q_{2}}-\Lambda_{q_{1}}\|+\displaystyle\frac{1}{\lambda}}\,\|\varphi\|_{H^{3}(\R^{n})}^{2}.
\end{equation}
holds true for any sufficiently large $\lambda>0.$
\end{Lemm}
\begin{Demo}{}
In view of Lemma \ref{Lem2.2} and using the fact that supp $\varphi\cap\Omega=\emptyset$, there exists  a geometrical optics solutions $u_{2,\lambda}$ to the equation
$$\para{\p_{t}^{2}-\Delta+q_{2}(t,x)}u_{2,\lambda}(t,x)=0 \,\,\,\mbox{in}\,\,Q, \,\,\,u_{2,\lambda|t=0}=\p_{t}u_{2,\lambda|t=0}=0\,\,\,\mbox{in}\,\,\Omega,$$
of the form
\begin{equation}\label{Eq3.15}
u_{2,\lambda}(t,x)=a(t,x)e^{i\lambda(x.\omega+t)}+R_{2,\lambda}(t,x),
\end{equation}
where $R_{2,\lambda}$ satisfies
$$\p_{t}R_{2,\lambda|t=0}=R_{2,\lambda|t=0}=0,\,\,\,\,R_{2,\lambda|\Sigma}=0.$$
and
\begin{equation}\label{Eq3.16}
\|R_{2,\lambda}\|_{L^{2}(Q)}\leq\displaystyle\frac{C}{\lambda}\,\|\varphi\|_{H^{3}(\R^{n})}.
\end{equation}
We denote by $u_{1}$, the solution of
$$
\left\{
  \begin{array}{ll}
    \para{\p_{t}^{2}-\Delta+q_{1}(t,x)}u_{1}(t,x)=0  & \mbox{in}\,\,Q, \\

    u_{1}(0,x)=\p_{t}u_{1}(0,x)=0 & \mbox{in}\,\,\Omega,\\

  u_{1}(t,x)=u_{2,\lambda}(t,x):=f_{\lambda}(t,x), & \mbox{on}\,\,\Sigma.
  \end{array}
\right.
$$
Putting $u(t,x)=u_{1}(t,x)-u_{2,\lambda}(t,x)$, we  get that
$$
\left\{
  \begin{array}{ll}
    \para{\p_{t}^{2}-\Delta+q_{1}(t,x)}u(t,x)=q(t,x)u_{2,\lambda}(t,x) & \mbox{in}\,\,Q, \\
    u(0,x)=\p_{t}u(0,x)=0 & \mbox{in}\,\,\Omega, \\
    u(t,x)=0 & \mbox{on}\,\,\Sigma.
  \end{array}
\right.
$$
Applying Lemma \ref{Lem2.2}, once more for $\lambda$  large enough and using the fact that supp $\varphi\pm T\omega\cap\Omega=\emptyset$, we may find a
geometrical optic solution $v_{\lambda}$ to the backward wave equation
$$\para{\p_{t}^{2}-\Delta+q_{1}(t,x)}v_{\lambda}(t,x)=0,\,\,\,\mbox{in}\,\, Q,\,\,\,\,v_{\lambda|t=T}=\p_{t}v_{\lambda|t=T}=0,\,\,\,\,\mbox{in}\,\,\Omega,$$
of the form
\begin{equation}\label{Eq3.17}
v_{\lambda}(t,x)=a(t,x)e^{-i\lambda(x.\omega+t)}+R_{1,\lambda}(t,x),
\end{equation}
where  $R_{1,\lambda}$ satisfies
$$\p_{t}R_{1,\lambda|t=T}=R_{1,\lambda|t=T}=0,\,\,\,\,R_{1,\lambda|\Sigma}=0,$$
and
\begin{equation}\label{Eq3.18}
\|R_{1,\lambda}\|_{L^{2}(Q)}\leq\frac{C}{\lambda}\,\|\varphi\|_{H^{3}(\R^{n})}.
\end{equation}
Consequently, by integrating by parts and using the Green's formula, we obtain
\begin{equation}\label{Eq3.19}
\begin{array}{lll}
\displaystyle\int_{Q}q(t,x)u_{2,\lambda}(t,x)v_{\lambda}(t,x)\,dx\,dt&=&\displaystyle\int_{Q}\para{\p_{t}^{2}-
\Delta+q_{1}(t,x)}u(t,x) v_{\lambda}(t,x)\,dx\,dt\\
&=&\displaystyle\int_{\Sigma} (\Lambda_{q_{2}}-\Lambda_{q_{1}})f_{\lambda}(t,x)v_{\lambda}(t,x)\,d\sigma\,dt,
\end{array}
\end{equation}
So, (\ref{Eq3.15}), (\ref{Eq3.17}) and  (\ref{Eq3.19}) yield
\begin{eqnarray}\label{Eq3.20}
\displaystyle\int_{Q}\!\!&q(t,x)&\!\!\!a^{2}(t,x)\,dx\,dt+\int_{Q}q(t,x)R_{1,\lambda}(t,x)R_{2,\lambda}(t,x)\,dx\,dt\cr
&&\!\!\!\!\!\!+\int_{Q}q(t,x)a(t,x)\para{R_{2,\lambda}(t,x)e^{-i\lambda(x.\omega+t)}+R_{1,\lambda}(t,x)e^{i\lambda(x.\omega+t)}}\,dx\,dt\cr
&=\displaystyle\int_{\Sigma}& \!\!\!(\Lambda_{q_{2}}-\Lambda_{q_{1}})f_{\lambda}(t,x)\,\,v_{\lambda}(t,x)\,\,d\sigma\,\,dt.
\end{eqnarray}
From (\ref{Eq3.20}), (\ref{Eq3.16}) and (\ref{Eq3.18}) it follows that
$$
|\int_{Q}q(t,x)\,a^{2}(t,x)\,dx\,dt|\leq\int_{\Sigma}|(\Lambda_{q_{2}}-\Lambda_{q_{1}})f_{\lambda}(t,x)\,v_{\lambda}(t,x)|\,d\sigma\,dt
+\frac{C}{\lambda}\,\|\varphi\|^{2}_{H^{3}(\R^{n})},
$$
where the constant $C>0$ does not depend on $\lambda$. Hence from the
Cauchy-Schwartz inequality and using the fact that $f_{\lambda}(t,x)=u_{2,\lambda}(t,x)$ on
$\Sigma$, we obtain
\begin{equation}\label{Eq3.21}
|\displaystyle\int_{Q}q(t,x)a^{2}(t,x)\,dx\,dt|\leq\|\Lambda_{q_{2}}-\Lambda_{q_{1}}\|\,\|u_{2,\lambda}\|_{H^{1}(\Sigma)}\,
\|v_{\lambda}\|_{L^{2}(\Sigma)}+\displaystyle\frac{C}{\lambda}\,\|\varphi\|^{2}_{H^{3}(\R^{n})},
\end{equation}
Further, as  $R_{i,\lambda|\Sigma}=0$, \,$\mbox{for}\,\,i=1,2,$ we
deduce from (\ref{Eq3.21}) that $$ |\int_{Q}q(t,x)\,a^{2}(t,x)\,dx\,dt|\leq C\para{
\|\Lambda_{q_{2}}-\Lambda_{q_{1}}\|\,\|u_{2,\lambda}-R_{2,\lambda}\|_{H^{2}(Q)}\,\|v_{\lambda}-R_{1,\lambda}\|_{H^{1}(Q)}+\frac{1}{\lambda}\,
\|\varphi\|^{2}_{H^{3}(\R^{n})}}.
$$
Bearing in mind that   $$\|v_{\lambda}-R_{1,\lambda}\|_{H^{1}(Q)}\leq
C\lambda\,\|\varphi\|_{H^{3}(\R^{n})},$$
$$\|u_{2,\lambda}-R_{2,\lambda}\|_{H^{2}(Q)}\leq C \lambda^{2}\,\|\varphi\|_{H^{3}(\R^{n})},$$
we end up getting that $$ |\int_{Q}q(t,x)\,a^{2}(t,x)\,dx\,dt|\leq
C\para{\lambda^{3}\|\Lambda_{q_{2}}-\Lambda_{q_{1}}\|+\frac{1}{\lambda}}\,\|\varphi\|^{2}_{H^{3}(\R^{n})}.
$$
Therefore by extending  $q(x,t)$ by zero outside $Q_{r}$ and recalling (\ref{Eq2.3}),
we find out that
$$
|\int_{0}^{T}\int_{\R^{n}}q(t,x-t\omega)\,\varphi^{2}(x)\,dx\,dt|\leq C\para{\lambda^{3}\|\Lambda_{q_{2}}-\Lambda_{q_{1}}\|+
\frac{1}{\lambda}}\|\varphi\|^{2}_{H^{3}(\R^{n})}.
$$
This completes the proof of the Lemma.
\end{Demo}
\subsection{X-ray transform}
The X-ray transform $R$ maps a function in $\R^{n+1}$ into the set of its
line integrals. More precisely, if $\omega\in \mathbb{S}^{n-1}$ and
$(t,x)\in\R^{n+1}$,
$$R(f)(\omega,x):=\int_{\R}f(t,x-t\omega) \,dt,$$
is the integral of $f$ over the lines $ \{(t,x-t\omega),\,\,\,t\in\R\}.$\\
\\
Using the above Lemma, we can estimate the X-ray transform of the differential of potentials as follows:
\begin{Lemm}\label{Lem3.2}
There exists a constant $C>0$, $\beta>0$, $\delta>0$, and $\lambda_{0}>0$
such that for all $\omega\in \mathbb{S}^{n-1},$ we have
$$|R(q)(\omega,y)|\leq C\,\para{\lambda^{\beta}\|\Lambda_{q_{2}}-\Lambda_{q_{1}}\|+\frac{1}{\lambda^{\delta}}},\,\,\,\,\mbox{a.e.}\,\,y\in\R^{n}.$$
for any $\lambda\geq\,\lambda_{0}.$
\end{Lemm}
\begin{Demo}{}
Let $\phi\in \!\mathcal{C}_{0}^{\infty}(\R^{n})$ be a positive function which
is supported in the unit ball $B(0,1)$ such that $\|\phi\|_{L^{2}(\R^{n})}=1.$ Define
$$\varphi_{\varepsilon}(x)=\varepsilon^{-n/2}\phi\para{\frac{x-y}{\varepsilon}}$$
where $y\in\mathscr{A}_{r}$. Then for sufficiently small $\varepsilon>0$ we can verify that
$$\mbox{supp}\,\varphi_{\varepsilon}\cap \Omega=\emptyset,\,\,\,\mbox{and}\,\,\,\,\mbox{supp}\,\varphi_{\varepsilon}\pm T\omega\cap \Omega=\emptyset.$$
And we have
$$\begin{array}{lll}
|\displaystyle\int_{0}^{T}\!\!\!\!&q&\!\!\!\!\!\!(t,y-t\omega)\,dt|=|\displaystyle\int_{0}^{T}\int_{\R^{n}}q(t,y-t\omega)\,
\varphi_{\varepsilon}^{2}(x)\,dx\,dt\,|\\
&\leq&|\displaystyle\int_{0}^{T}\int_{\R^{n}}q(t,x-t\omega)\,\varphi^{2}_{\varepsilon}(x)\,dx\,dt|+|\displaystyle\int_{0}^{T}
\int_{\R^{n}}\displaystyle\left(q(t,y-t\omega)
-q(t,x-t\omega)\right)\,\varphi_{\varepsilon}^{2}(x)\,dx\,dt\,|.
\end{array}$$
Since $\|q\|_{C^{1}(Q)}\leq M$, we have
$$|q(t,y-t\omega)-q(t,x-t\omega)|\leq C \,|x-y|.$$
Applying Lemma \ref{Lem3.1} with $\varphi=\varphi_{\varepsilon}$, we obtain
\begin{equation}\label{Eq3.22}
|\int_{0}^{T}q(t,y-t\omega)\,dt|\leq C\para{\lambda^{3}\|\Lambda_{q_{2}}-\Lambda_{q_{1}}\|+\frac{1}{\lambda}}\|\varphi_{\varepsilon}\|^{2}_{H^{3}(\R^{n})}
+ C \int_{\R^{n}}|x-y|\,\varphi_{\varepsilon}^{2}(x)\,dx.
\end{equation}
On the other hand, we have
$$\|\varphi_{\varepsilon}\|_{H^{3}(\R^{n})}\leq C\,\varepsilon^{-3},\,\,\,\,\,\int_{\R^{n}}|x-y|\,\varphi_{\varepsilon}^{2}(x)\,dx\leq C \varepsilon.$$
Thus, from (\ref{Eq3.22}), we obtain
$$|\int_{0}^{T}q(t,y-t\omega)\,dt|\leq C\para{\lambda^{3}\|\Lambda_{q_{2}}-\Lambda_{q_{1}}\|+\frac{1}{\lambda}} \varepsilon^{-6}+C\varepsilon. $$
We select $\varepsilon$ such that
$$\varepsilon=\frac{\varepsilon^{-6}}{\lambda}.$$
Then there exist constants $\delta>0$ and $\beta>0$ such that
$$|\int_{0}^{T}q(t,y-t\omega)\,dt|\leq\,C\,\para{\lambda^{\beta}\|\Lambda_{q_{2}}-\Lambda_{q_{1}}\|+\frac{1}{\lambda^{\delta}}}.$$
Using the fact that \, $q=0$, outside $Q_{r}$, we get
\begin{equation}\label{Eq3.23}
|\int_{\R}q(t,y-t\omega)\,dt|\leq\,C\,\para{\lambda^{\beta}\|\Lambda_{q_{2}}-\Lambda_{q_{1}}\|+\frac{1}{\lambda^{\delta}}},\,\,\,\,\,\,
\mbox{a.e.}\,y\in\mathscr{A}_{r},\,\,\omega\in\mathbb{S}^{n-1}.
\end{equation}
On the other hand, if  $|y|\leq \displaystyle\frac{r}{2}$, then
\begin{equation}\label{**}
q(t,y-t\omega)=0\,\,\forall\,t\in\R.
\end{equation}
Indeed, we have
\begin{equation}\label{*}
|y-t\omega|\geq|t|-|y|\geq t-\displaystyle\frac{r}{2}.
\end{equation}
So that, if $t>\displaystyle\frac{r}{2}$, from (\ref{*}), we have $(t,y-t\omega)\notin \mathscr{C}_{r}^{+}$. And if $t\leq \displaystyle\frac{r}{2}$, we have also $(t,y-t\omega)\notin \mathscr{C}_{r}^{+}$.
Consequently,
$$(t,y-t\omega)\notin \mathscr{C}_{r}^{+}\supset Q_{*},\,\,\,\,\,\,\mbox{for\,all} \,\,t\in\R.$$
Using the fact that $q=q_{2}-q_{1}=0$ outside $Q_{*}$, we deduce (\ref{**}). Therefore,
$$\int_{\R}q(t,y-t\omega)\,dt=0,\,\,\,\,\mbox{a.e.}\,\,y\in B(0,\frac{r}{2}).$$
By a similar way, we prove that in the case where $|y|\geq T-\frac{r}{2}$, we have
$$(t,y-t\omega)\notin\mathscr{C}_{r}^{-}\supset Q_{*},\,\,\,\,\,\mbox{for\,all}\,\,t\in\R.$$
Then we conclude that
\begin{equation}\label{Eq3.24}
\int_{\R}q(t,y-t\omega)\,dt=0,\,\,\,\,\mbox{a.e.}\,\,y\notin\mathscr{A}_{r},\,\,\,\omega\in \mathbb{S}^{n-1}.
\end{equation}
Consequently, by (\ref{Eq3.23}) and (\ref{Eq3.24}), one gets
$$
|R(q)(\omega,y)|=|\int_{\R}q(t,y-t\omega)\,dt|\leq\, C\,\para{\lambda^{\beta}\|\Lambda_{q_{2}}-\Lambda_{q_{1}}\|+\frac{1}{\lambda^{\delta}}},\,\,\,\mbox{a.e}.\,\,y\in\R^{n},\,\,\,\omega\in\mathbb{S}^{n-1}.
$$
This completes the proof of the Lemma.
\end{Demo}
Let now $$E=\{(\tau,\xi)\in\R\times\R^{n},\,\,\,|\tau|\leq|\xi|\},$$ and let the Fourier transform of $q\in L^{1}(\R^{n+1})$
$$\widehat{q}(\tau,\xi)=\int_{\R}\int_{\R^{n}}q(x,t)e^{-ix.\xi}e^{-it\tau}\,dx\,dt.$$
Our goal now is to prove the following
\begin{Lemm}\label{Lem3.3}
There exist  constants $C>0$, $\beta>0,$ $\delta>0$ and $\lambda_{0}>0$ such
that the following estimate holds
$$|\widehat{q}(\tau,\xi)|\leq C\,\para{\lambda^{\beta}\|\Lambda_{q_{2}}-\Lambda_{q_{1}}\|+\frac{1}{\lambda^{\delta}}},$$
 for any $(\tau,\xi)\in E$ and $\lambda\geq \lambda_{0}$.
\end{Lemm}
\begin{Demo}{}
Let $(\tau,\xi)\in E$ and $\zeta\in\mathbb{S}^{n-1}$ such that $\xi.\zeta=0$. By defining
$$\omega=\frac{\tau}{|\xi|^{2}}.\xi+\sqrt{1-\frac{\tau^{2}}{|\xi|^{2}}}.\zeta,$$
we have $\omega\in\mathbb{S}^{n-1}$ and $\omega.\xi=\tau.$\\
\\
By the change of variable $x=y-t\omega$ we have for all $\xi\in\R^{n},\,\,\omega\in\mathbb{S}^{n-1}$
$$\begin{array}{lll}
\displaystyle\int_{\R^{n}}R(q)(\omega,y)\,e^{-iy.\xi}\,dy&=&\displaystyle\int_{\R^{n}}\displaystyle\para{\int_{\R}q(t,y-t\omega)\,dt}\,e^{-iy.\xi}\,dy\\
&=&\displaystyle\int_{\R}\int_{\R^{n}}q(t,x)\,e^{-ix.\xi}e^{-it(\omega.\xi)}\,dx\,dt\\
&=&\widehat{q}(\omega.\xi,\xi)\\
&=&\widehat{q}(\tau,\xi).
\end{array}$$
Denote $(\tau,\xi)=(\omega.\xi,\xi)\in E$. Since
$\mbox{supp}\,q(t,.)\subset\Omega\subset B(0,\frac{r}{2})$, then we have
 $$\int_{\R^{n}\cap B(0,\frac{r}{2}+T)}R(q)(\omega,y)\,e^{-iy.\xi}\,dy=\widehat{q}(\tau,\xi).$$
In terms of Lemma \ref{Lem3.2}, the proof is completed.
\end{Demo}
\subsection{Stability estimate}
We are now in position to complete the proof of Theorem \ref{Thm1}.
For $\rho >0$ and  $\gamma\in
(\mathbb{N}\cup\{0\})^{n+1}$, we denote
 $$|\gamma|=\gamma_{1}+...+\gamma_{n+1},\,\,\,\,\,\,\,\,\,B(0,\rho)=\{x\in\R^{n+1},\,\,|x|<\rho\}.$$
We consider the following Lemma
\begin{Lemm}(see \cite{[C4]})\label{Lem3.4}
Let $O$ be an open set of $B(0,1)$, and $F$ an analytic function in $B(0,2),$
satisfying the following property: there exist constant $M,\eta>0$ such that
$$\|\p^{\gamma}F\|_{L^{\infty}(B(0,2))}\leq \frac{M|\gamma|!}{\eta^{|\gamma|}},\,\,\,\,\forall\,\gamma\in(\mathbb{N}\cup\{0\})^{n+1}.$$
Then,
$$\|F\|_{L^{\infty}(B(0,1))}\leq (2M)^{1-\mu}\|F\|_{L^{\infty}(O)}^{\mu}.$$
where $\mu\in(0,1)$ depends on $n$, $\eta$ and $|O|$.
\end{Lemm}

The Lemma is conditional stability for the analytic continuation, and see Lavrent'ev, Romanov and Shishat$\cdot$sKii. \cite{[M.M]} for classical results.  For fixed $\alpha>0$, let us set
$$F_{\alpha}(\tau,\xi)=\widehat {q}(\alpha(\tau,\xi))\,\,\, \mbox{for}\,\,\, (\tau,\xi)\in\R^{n+1}.$$
It is easily seen that $F_{\alpha}$ is analytic and we have
\begin{eqnarray}\label{Eq3.25}
|\p^{\gamma} F_ {\alpha}(\tau,\xi)|=|\p^{\gamma}
\widehat{q}(\alpha(\tau,\xi))|&=&|\p^{\gamma}\displaystyle\int_{\R^{n+1}}q(t,x)\,e^{-i\alpha
(t,x).(\tau,\xi)}\,dx\,dt|\cr
&=&|\displaystyle\int_{\R^{n+1}}q(t,x)(-i)^{|\gamma|}\alpha^{|\gamma|}(t,x)^{\gamma}e^{-i\alpha(t,x).(\tau,\xi)}\,dx\,dt|.
\end{eqnarray}
Therefore, from (\ref{Eq3.25}) one gets
$$\begin{array}{lll}
|\p^{\gamma} F_ {\alpha}(\tau,\xi)|\leq \displaystyle\int_{\R^{n+1}}|q(t,x)| \alpha^{|\gamma|}(|x|^{2}+t^{2})^{\frac{|\gamma|}{2}}\,dx\,dt\leq\|q\|_{L^{1}(Q_{*})}\,\,\alpha^{|\gamma|}\,\,(2T^{2})^{\frac{|\gamma|}{2}}
\leq C \,\,\displaystyle\frac{|\gamma|!}{(T^{-1})^{|\gamma|}}\,\,e^{\alpha}.
\end{array}
$$
Then, applying Lemma \ref{Lem3.4} in the set $O=\mathring{E}\cap B(0,1)$ with
$M=Ce^{\alpha}$ and $\eta=T^{-1},$ we can take a constant $\mu\in(0,1)$ such
that
$$|F_{\alpha}(\tau,\xi)|=|\widehat{q}(\alpha(\tau,\xi))|\leq C e^{\alpha(1-\mu)}\|F_{\alpha}\|_{L^{\infty}(O)}^{\mu},\,\,\,\,\,\,\,\,\,(\tau,\xi)\in B(0,1). $$
Hence, by using the fact that $\alpha\,\mathring{E}=\{\alpha(
\tau,\xi),\,(\tau,\xi)\in\mathring{E}\}=\mathring{E}$, we get for $(\tau,\xi)\in B(0,\alpha)$
\begin{eqnarray}\label{Eq3.26}
|\widehat{q}(\tau,\xi)|=|F_{\alpha}(\alpha^{-1}(\tau,\xi)|&\leq& C e^{\alpha(1-\mu)}\,\|F_{\alpha}\|_{L^{\infty}(O)}^{\mu}\cr
&\leq& C e^{\alpha(1-\mu)} \|\,\widehat{q}\,\|^{\mu}_{L^{\infty}(B(0,\alpha)\cap\mathring{E})}\cr
&\leq&C e^{\alpha (1-\mu)}\|\,\widehat{q}\,\|_{L^{\infty}(\mathring{E})}^{\mu}.
\end{eqnarray}
On the other hand we have
$$\begin{array}{lll}
\|q\|_{H^{-1}(\R^{n+1})}^{2/\mu}\!\!\!&=\!\!\!&\displaystyle\para{\displaystyle\int_{|(\tau,\xi)|<\alpha}\!\!\!\!(1+|(\tau,\xi)|^{2})^{-1}
|\widehat{q}(\tau,\xi)|^{2}\,d\tau d\xi
+\!\!\displaystyle\int_{|(\tau,\xi)|\geq\alpha}\!\!\!\!\!(1+|(\tau,\xi)|^{2})^{-1}|\widehat{q}(\tau,\xi)|^{2}\,d\tau d\xi\,}^{1/\mu}\\
&\leq&  C\para{\alpha^{n+1}\,\,\,\|\widehat{q}\|^{2}_{L^{\infty}(B(0,\alpha))}+\,\alpha^{-2}\,\,\|q\|^{2}_{L^{2}(\R^{n+1})}}^{1/\mu}.\\
\end{array}$$
From (\ref{Eq3.26}) and applying Lemma \ref{Lem3.3}, we obtain
$$\begin{array}{lll}
\|q\|^{2/\mu}_{H^{-1}(\R^{n+1})}&\leq&\,C\displaystyle\para{\alpha^{{n+1}}\,e^{2\alpha(1-\mu)}\,(\lambda^{\beta}\|\Lambda_{q_{2}}-\Lambda_{q_{1}}\|
+\frac{1}{\lambda^{\delta}})^{
2\mu}+\alpha^{-2}}^{1/\mu}\\
&\leq&C\displaystyle\para{\alpha^{\frac{n+1}{\mu}} \,e^{\frac{2\alpha(1-\mu)}{\mu}} \lambda^{2\beta}\|\Lambda_{q_{2}}-\Lambda_{q_{1}}\|^{2}+\alpha^{\frac{n+1}{\mu}}\,e^{\frac{2\alpha(1-\mu)}{\mu}}\,\lambda^{-2\delta}+\alpha^{-2/\mu}}.
\end{array}$$
Let $\alpha_{0}>0$ be sufficiently large and $\alpha>\alpha_{0}$. Set
$$\lambda=\alpha^{\frac{n+3}{2\mu\delta}}\,e^{\frac{\alpha(1-\mu)}{\mu\delta}}.$$
By $\alpha>\alpha_{0},$ we can assume that $\lambda>\lambda_{0}$, and we have
$$\alpha^{\frac{n+1}{\mu}}\,e^{\frac{2\alpha(1-\mu)}{\mu}}\,\lambda^{-2\delta}=\alpha^{-2/\mu}.$$
Then $$\begin{array}{lll} \|q\|^{2/\mu}_{H^{-1}(\R^{n+1})}&\leq&
C\para{\alpha^{\frac{\delta(n+1)+\beta(n+3)}{\delta\mu}}\,e^{\frac{2\alpha(\delta+\beta)(1-\mu)}{\delta\mu}}\|\Lambda_{q_{2}}-\Lambda_{q_{1}}\|^{2}
+\alpha^{-2/\mu}}\\
&\leq&
C\,\displaystyle\para{e^{N\alpha}\|\Lambda_{q_{2}}-\Lambda_{q_{1}}\|^{2}+\alpha^{-2/\mu}},
\end{array}$$
where $N$ depends on $\delta,\,\beta,\,n,$ and $\mu$. In order to minimize
the right hand-side  with respect to  $\alpha$, we set
 $$\alpha=\frac{1}{N}|\,\log\|\Lambda_{q_{2}}-\Lambda_{q_{1}}\|\,|,$$
  where we assume that
$$0<\|\Lambda_{q_{2}}-\Lambda_{q_{1}}\|<c.$$
  It follows that
$$\begin{array}{lll}\|q\|_{H^{-1}(Q_{*})}\leq \|q\|_{H^{-1}(\R^{n+1})}&\leq&
C\para{\|\Lambda_{q_{2}}-\Lambda_{q_{1}}\|+|\,\log\|\Lambda_{q_{2}}-\Lambda_{q_{1}}\|\,|^{-2/\mu}}^{\mu/2}\\
&\leq&C\para{\|\Lambda_{q_{2}}-\Lambda_{q_{1}}\|^{\mu/2}+|\log\|\Lambda_{q_{2}}-\Lambda_{q_{1}}\||^{-1}}.\end{array}$$
The estimate (\ref{Eq1.2}), is now an easy consequence of  the Sobolev embedding
theorem and an interpolation inequality.  Let $\delta'>0$ such that
$s=n/2+2\delta'$. Then, we have 
$$
\begin{array}{lll}
\|q\|_{L^{\infty}(Q_{*})}&\leq&C\|q\|_{H^{s}(Q_{*})}\\
&\leq& C \,\|q\|_{H^{-1}(Q_{*})}^{1-\beta}\,\|q\|_{H^{s+1}(Q_{*})}^{\beta}\\
&\leq&C\,\|q\|_{H^{-1}(Q_{*})}^{1-\beta},
\end{array}$$
for some $\beta\!\in\!(0,1)$. 
 Then the proof of Theorem \ref{Thm1} is completed.
\section{Proof of Theorem \ref{Thm2}}\label{Sec4}
This section is devoted to the proof of Theorem \ref{Thm2}. We will extend the stability estimate (\ref{Eq1.2}) given in Theorem \ref{Thm1}, to an estimate in a larger region $Q_{\sharp}\supset Q{*}$. Differently  to Theorem \ref{Thm1}, here the observations are given by the boundary operator $\mathscr{R}_{q}$ introduced in Subsection \ref{Sub1.2}. We need to consider geometric optics solutions similar to the one used in the previous section, but this time, we will only assume that $\mbox{supp}\,\varphi\cap\Omega=\emptyset.$
(We don't need to assume that supp $\varphi\pm T\omega\cap\Omega=\emptyset$). Let's first recall the definition of the operator $\mathscr{R}_{q}$:
$$\begin{array}{ccc}
\mathscr{R}_{q}:H^{1}(\Sigma)&\longrightarrow& L^{2}(\Sigma)\times H^{1}(\Omega)\times L^{2}(\Omega).\\
f&\longmapsto&(\p_{\nu}u,\,u(T,.),\,\p_{t}u(T,.)).
\end{array}$$
We denote by
$$\mathscr{R}_{q_{j}}^{1}(f)=\p_{\nu}^{\,}u_{j},\,\,\,\,\mathscr{R}_{q_{j}}^{2}(f)=u_{j,}(T,.),\,\,\,
\mathscr{R}_{q_{j}}^{3}(f)=\p_{t}u_{j}(T,.),\,\,\,\,\mbox{for}\,j=1,2.$$
\begin{Lemm}\label{Lem4.1} Let $q_{1},q_{2}\in \mathcal{A}^{\sharp}(q_{0},M)$, $\varphi\in \mathcal{C}^{\infty}_{0}(\R^{n})$, such that supp\,$\varphi\cap\Omega=\emptyset$,  and put $q=\para{q_{2}-q_{1}}$. Then, there exists $C>0,$ such that for any $\omega\in \mathbb{S}^{n-1}$
the following estimate
\begin{equation}\label{Eq4.27}
|\displaystyle\int_{0}^{T}\int_{\R^{n}}q(t,x-t\omega)\,\varphi^{2}(x)\,dx\,dt|\leq C\,\para{\lambda^{3}\|\mathscr{R}_{q_{2}}-\mathscr{R}_{q_{1}}\|+\displaystyle\frac{1}{\lambda}}\,\|\varphi\|_{H^{3}(\R^{n})}^{2}
\end{equation}
holds true for any sufficiently large $\lambda>0.$
\end{Lemm}
\begin{Demo}{}
In view of Lemma \ref{Lem2.2} and using the fact that supp $\varphi\cap\Omega=\emptyset$, there exists  a geometrical optics solutions $u_{2,\lambda}$ to the equation
$$\para{\p_{t}^{2}-\Delta+q_{2}(t,x)}u_{2,\lambda}(t,x)=0 \,\,\,\mbox{in}\,\,Q, \,\,\,u_{2,\lambda|t=0}=\p_{t}u_{2,\lambda|t=0}=0\,\,\,\mbox{in}\,\,\Omega,$$
of the form
\begin{equation}\label{Eq4.28}
u_{2,\lambda}(t,x)=a(t,x)e^{i\lambda(x.\omega+t)}+R_{2,\lambda}(t,x),
\end{equation}
where $R_{2,\lambda}$ satisfies
$$\p_{t}R_{2,\lambda|t=0}=R_{2,\lambda|t=0}=0,\,\,\,\,R_{2,\lambda|\Sigma}=0,$$
and
\begin{equation}\label{Eq4.29}
\|R_{2,\lambda}\|_{L^{2}(Q)}\leq\displaystyle\frac{C}{\lambda}\,\|\varphi\|_{H^{3}(\R^{n})}.
\end{equation}
We denote by $u_{1,\lambda}$, the solution of
$$
\left\{
  \begin{array}{ll}
    \para{\p_{t}^{2}-\Delta+q_{1}(t,x)}u_{1,\lambda}(t,x)=0  & \mbox{in}\,\,Q, \\

    u_{1,\lambda}(0,x)=\p_{t}u_{1,\lambda}(0,x)=0 & \mbox{in}\,\,\Omega,\\

  u_{1,\lambda}(t,x)=u_{2,\lambda}(t,x):=f_{\lambda}(t,x), & \mbox{on}\,\,\Sigma.
  \end{array}
\right.
$$
Putting $u_{\lambda}(t,x)=u_{1,\lambda}(t,x)-u_{2,\lambda}(t,x)$, we  get that
$$
\left\{
  \begin{array}{ll}
    \para{\p_{t}^{2}-\Delta+q_{1}(t,x)}u_{\lambda}(t,x)=q(t,x)u_{2,\lambda}(t,x) & \mbox{in}\,\,Q \\
    u_{\lambda}(0,x)=\p_{t}u_{\lambda}(0,x)=0 & \mbox{in}\,\,\Omega \\
    u_{\lambda}(t,x)=0 & \mbox{on}\,\,\Sigma.
  \end{array}
\right.
$$
Applying Lemma \ref{Lem2.2}, once more for $\lambda$  large enough,  we may find a
geometrical optic solution $v_{\lambda}$ to the backward wave equation
$$\para{\p_{t}^{2}-\Delta+q_{1}(t,x)}v_{\lambda}(t,x)=0,\,\,\,\mbox{in}\,\, Q,$$
of the form
\begin{equation}\label{Eq4.30}
v_{\lambda}(t,x)=a(t,x)e^{-i\lambda(x.\omega+t)}+R_{1,\lambda}(t,x),
\end{equation}
where  $R_{1,\lambda}$ satisfies
$$\p_{t}R_{1,\lambda|t=T}=R_{1,\lambda|t=T}=0,\,\,\,\,R_{1,\lambda|\Sigma}=0,$$
and
\begin{equation}\label{Eq4.31}
\|R_{1,\lambda}\|_{L^{2}(Q)}\leq\frac{C}{\lambda}\,\|\varphi\|_{H^{3}(\R^{n})}.
\end{equation}
Consequently, by integrating by parts and using the Green's formula we obtain
\begin{eqnarray}\label{Eq4.32}
\displaystyle\int_{Q}q(t,x)u_{2,\lambda}(t,x)v_{\lambda}(t,x)\,dx\,dt&=&
\displaystyle\int_{\Sigma} (\mathscr{R}^{1}_{q_{2}}-\mathscr{R}^{1}_{q_{1}})(f_{\lambda})v_{\lambda}(t,x)\,d\sigma\,dt\cr
&&+\displaystyle\int_{\Omega}\displaystyle\para{\mathscr{R}_{q_{2}}^{2}-\mathscr{R}_{q_{1}}^{2}}(f_{\lambda})\,\p_{t}v_{\lambda}(T,.)\,dx\cr
&&-\displaystyle\int_{\Omega}\para{\mathscr{R}_{q_{2}}^{3}-\mathscr{R}_{q_{1}}^{3}}(f_{\lambda})\,v_{\lambda}(T,.)\,dx,
\end{eqnarray}
Then, by replacing $u_{2,\lambda}$ and $v_{\lambda}$ by their expressions in the left hand side of (\ref{Eq4.32}) and using (\ref{Eq4.29}) and (\ref{Eq4.31}), then from Cauchy-Schwartz inequality, one gets the following estimate
$$\begin{array}{lll}
|\displaystyle\int_{Q}q(t,x)\,a^{2}(t,x)\,dx\,dt|&\leq&\|(\mathscr{R}_{q_{2}}^{1}-\mathscr{R}_{q_{1}}^{1})(f_{\lambda})\|_{L^{2}(\Sigma)}
\|v_{\lambda}\|_{L^{2}(\Sigma)}+
\displaystyle\frac{C}{\lambda}\|\varphi\|^{2}_{H^{3}(\R^{n})}\\
\\
&&+\|(\mathscr{R}_{q_{2}}^{2}-\mathscr{R}^{2}_{q_{1}})(f_{\lambda})\|_{L^{2}(\Omega)}\|\p_{t}v_{\lambda}(T,.)\|_{L^{2}(\Omega)}\\
\\
&&+\|(\mathscr{R}_{q_{2}}^{3}-\mathscr{R}_{q_{1}}^{3})(f_{\lambda})\|_{L^{2}(\Omega)}\|v_{\lambda}(T,.)\|_{L^{2}(\Omega)}.
\end{array}$$
Then we obtain,
\begin{eqnarray}\label{Eq4.33}
\!\!|\!\displaystyle\int_{Q}\!q(t,x)\,a^{2}(t,x)\,dx\,dt|\!\!\!&\leq&\!\!\!\!\displaystyle\left( \|(\mathscr{R}_{q_{2}}^{1}\!-\mathscr{R}_{q_{1}}^{1})(f_{\lambda})\|^{2}_{L^{2}(\Sigma)}+\|(\mathscr{R}_{q_{2}}^{2}\!-\mathscr{R}_{q_{1}}^{2})
(f_{\lambda})\|^{2}_{H^{1}(\Omega)}
+\|(\mathscr{R}_{q_{2}}^{3}\!-\mathscr{R}_{q_{1}}^{3})(f_{\lambda})\|^{2}_{L^{2}(\Omega)}\right)^{\frac{1}{2}}\cr
&&\displaystyle\left( \|v_{\lambda}\|^{2}_{L^{2}(\Sigma)}+\|v_{\lambda}(T,.)\|_{L^{2}(\Omega)}^{2}+\|\p_{t}v_{\lambda}(T,.)\|^{2}_{L^{2}(\Omega)}  \right)^{\frac{1}{2}}+\displaystyle\frac{C}{\lambda}\|\varphi\|^{2}_{H^{3}(\R^{n})}.
\end{eqnarray}
Setting
$$\phi_{\lambda}=(v_{\lambda|_{\Sigma}},\,v_{\lambda}(T,.),\,\p_{t}v_{\lambda}(T,.))$$
Then, from (\ref{Eq4.33}), we get
$$|\!\int_{Q}\!q(t,x) a^{2}(t,x)\,dx\,dt|\leq\|(\mathscr{R}_{q_{2}}-\mathscr{R}_{q_{1}})(f_{\lambda})\|_{L^{2}(\Sigma)\times H^{1}(\Omega)\times L^{2}(\Omega)}\,\,\|\phi_{\lambda}\|_{L^{2}(\Sigma)\times L^{2}(\Omega)\times L^{2}(\Omega)}+\frac{C}{\lambda}\|\varphi\|^{2}_{H^{3}(\R^{n})}$$
and using the fact that $f_{\lambda}(t,x)=u_{2,\lambda}(t,x)$ on
$\Sigma$, we obtain
$$\begin{array}{lll}
|\displaystyle\int_{Q}q(t,x)a^{2}(t,x)\,dx\,dt|
\leq\|\mathscr{R}_{q_{2}}-\mathscr{R}_{q_{1}}\|\,\|u_{2,\lambda}\|_{H^{1}(\Sigma)}\,\|\phi_{\lambda}\|_{L^{2}(\Sigma)\times L^{2}(\Omega)\times L^{2}(\Omega)}+\displaystyle\frac{C}{\lambda}\,\|\varphi\|^{2}_{H^{3}(\R^{n})},
\end{array}$$
Further, as  $R_{i,\lambda|\Sigma}=0$, $\mbox{for}\,\,i=1,2,$ we
deduce that $$ |\int_{Q}q(t,x)\,a^{2}(t,x)\,dx\,dt|\leq C\para{
\|\mathscr{R}_{q_{2}}-\mathscr{R}_{q_{1}}\|\,\|u_{2,\lambda}-R_{2,\lambda}\|_{H^{2}(Q)}\,\|\phi_{1,\lambda}\|_{H^{1}(Q)\times L^{2}(\Omega)\times L^{2}(\Omega)}+\frac{1}{\lambda}\,\|\varphi\|^{2}_{H^{3}(\R^{n})}},
$$
where $$\phi_{1,\lambda}=\displaystyle\left(v_{\lambda}-R_{1,\lambda}, \,v_{\lambda}(T,.),\,\p_{t}v_{\lambda}(T,.)\right). $$
Using the fact that $R_{1,\lambda}(T,.)=\p_{t}R_{1,\lambda}(T,.)=0$\, on $\Omega$, we have
  $$\|u_{2}-R_{2}\|_{H^{2}(Q)}\leq C \lambda^{2}\,\|\varphi\|_{H^{3}(\R^{n})},$$
and  $$\begin{array}{lll}
   \|\phi_{1,\lambda}\|_{H^{1}(Q)\times L^{2}(\Omega)\times L^{2}(\Omega)}&\leq&\|v_{\lambda}-R_{1,\lambda}\|_{H^{1}(Q)}+\|v_{\lambda|t=T}\|_{L^{2}(\Omega)}+\|\p_{t}v_{\lambda|t=T}\|_{L^{2}(\Omega)}\\
   &\leq& C\lambda\,\|\varphi\|_{H^{3}(\R^{n})},
   \end{array}$$
Therefore by extending  $q(t,x)$ by zero outside $Q_{r}$ and recalling (\ref{Eq2.3}), we find out that
$$
|\int_{0}^{T}\int_{\R^{n}}q(t,x-t\omega)\,\varphi^{2}(x)\,dx\,dt|\leq C\para{\lambda^{3}\|\mathscr{R}_{q_{2}}-\mathscr{R}_{q_{1}}\|+
\frac{1}{\lambda}}\|\varphi\|^{2}_{H^{3}(\R^{n})}.
$$
This completes the proof of the Lemma.
\end{Demo}
Let's move now to prove the following Lemma
\begin{Lemm}\label{Lem4.2}
There exists a constant $C>0$, $\beta>0$, $\delta>0$, and $\lambda_{0}>0$
such that for all $\omega\in \mathbb{S}^{n-1},$ we have
$$|R(q)(\omega,y)|\leq C\,\para{\lambda^{\beta}\|\mathscr{R}_{q_{2}}-\mathscr{R}_{q_{1}}\|+\frac{1}{\lambda^{\delta}}},\,\,\,\,\mbox{a.e.}\,\,y\in\R^{n}.$$
for any $\lambda\geq\,\lambda_{0}.$
\end{Lemm}
\begin{Demo}{}
We consider $(\varphi_{\varepsilon})_{\varepsilon}$ defined in the proof of Lemma \ref{Lem3.2}. We only assume that  $y\notin\Omega$, then for  sufficiently small $\varepsilon>0$, we can verify that supp $\varphi_{\varepsilon}\cap\Omega=\emptyset$. Taking in acount this last remark, using Lemma \ref{Lem4.1} and repeating the arguments used in Lemma \ref{Lem3.2}, we obtain this estimate
\begin{equation}\label{Eq4.34}
|\int_{\R}q(t,y-t\omega)\,dt|\leq\,C\,\para{\lambda^{\beta}\|\mathscr{R}_{q_{2}}-\mathscr{R}_{q_{1}}\|+\frac{1}{\lambda^{\delta}}},\,\,\,\,\,\,
\mbox{a.e.}\,y
\notin B(0,\frac{r}{2}).
\end{equation}
On the other hand, if  $y\in B(0,\frac{r}{2})$, then we have 
\begin{equation}\label{Es2}
q(t,y-t\omega)=0,\,\,\,\forall\,t\in\R.
\end{equation}
Indeed, we have
\begin{equation}\label{Estim}
|y-t\omega|\geq|t|-|y|\geq t-\displaystyle\frac{r}{2}.
\end{equation}
So that, from (\ref{Estim}), we deduce that for all $t>\frac{r}{2}$ we have  $(t,y-t\omega)\notin \mathscr{C}_{r}^{+}.$ And if $t\leq \frac{r}{2}$, we have also that $(t,y-t\omega)\notin \mathscr{C}^{+}_{r}$. We recall that $Q_{\sharp}=Q\cap\mathscr{C}^{+}_{r}$.
 Consequently, we have
  $$(t,y-t\omega)\notin Q_{\sharp},\,\,\,\mbox{for\, all}\,\, t\in\R.$$
Then, using the fact that $q=q_{2}-q_{1}=0$ outside $Q_{\sharp}$,  we obtain (\ref{Es2}). Therefore
\begin{equation}\label{Eq4.35}
\int_{\R}q(t,y-t\omega)\,dt=0,\,\,\,\,\mbox{a.e.}\,\,y\in B(0,\frac{r}{2}).
\end{equation}
In light of (\ref{Eq4.34}) and (\ref{Eq4.35}), the proof of Lemma \ref{Lem4.2} is completed.
\end{Demo}

Using the above result and in the same way as in Section \ref{Sec3}, we complete the proof of Theorem \ref{Thm2}.
\section{Proof of Theorem \ref{Thm3}}\label{Sec5}
In this section we deal with the same problem treated in Section \ref{Sec3} and \ref{Sec4}, except our data will be the response of the medium for all possible initial data. As usual, we will prove Theorem \ref{Thm3} using geometric optics solutions constructed in Section \ref{Sec2} and X-ray transform. Let's first recall the definition of the operator $\mathcal{I}_{q}$:
$$\begin{array}{ccc}
\mathcal{I}_{q}: H^{1}(\Sigma)\times H^{1}(\Omega)\times L^{2}(\Omega)&\longrightarrow& L^{2}(\Sigma)\times H^{1}(\Omega)\times L^{2}(\Omega).\\
\psi=(f,u_{0},u_{1})&\longmapsto&(\p_{\nu}u,\,u(T,.),\,\p_{t}u(T,.)).
\end{array}$$
We denote by
$$\mathcal{I}_{q_{j}}^{1}(\psi)=\p_{\nu}^{\,}u_{j},\,\,\,\,\mathcal{I}_{q_{j}}^{2}(\psi)=u_{j}(T,.),\,\,\,\,
\mathcal{I}_{q_{j}}^{3}(\psi)=\p_{t}u_{j}(T,.),\,\,\,\,\mbox{for}\,j=1,2.$$
\begin{Lemm}\label{Lem5.1}
Let $q_{1},\,q_{2}\in \mathcal{C}^{1}(\overline{Q})$, and put $q=(q_{2}-q_{1})$. There exists $C>0$, $\beta>0,\,\,\,\delta>0$ and $\lambda_{0}>0$ such that for any $\omega\in\mathbb{S}^{n-1}$ we have  the following estimate
$$|R(q)(\omega,y)|\leq C\para{\lambda^{\beta}\|\mathcal{I}_{q_{2}}-\mathcal{I}_{q_{1}}\|+\frac{1}{\lambda^{\delta}}},\,\,\,\,\,\mbox{a.e}.\,y\in\R^{n}.$$
for any $\lambda\geq \lambda_{0}$.
\end{Lemm}
\begin{Demo}{}
Let $\varphi\in\mathcal{C}_{0}^{\infty}(\R^{n}).$ For $\lambda$ sufficiently large, Lemma \ref{Lem2.2} guarantees the existence of the geometrical optics solution $u_{2,\lambda}$ to
$$(\p_{t}^{2}-\Delta+q_{2}(t,x))u_{2,\lambda}(t,x)=0,\,\,\,\,\,\,\,\,\,\,\mbox{in}\,\,Q,$$
of the form
\begin{equation}\label{Eq5.36}
u_{2,\lambda}(t,x)=a(t,x)e^{i\lambda(x.\omega+t)}+R_{2,\lambda}(t,x)
\end{equation}
where $R_{2,\lambda}$ satisfies
$$\p_{t}R_{2,\lambda|t=0}=R_{2,\lambda|t=0}=0,\,\,\,\,R_{2,\lambda|\Sigma}=0,$$
and
\begin{equation}\label{Eq5.37}
\|R_{2,\lambda}\|_{L^{2}(Q)}\leq\frac{C}{\lambda}\|\varphi\|_{H^{3}(\R^{n})}.
\end{equation}
We denote $u_{1,\lambda}$ the solution of
$$
\left\{
  \begin{array}{ll}
    \para{\p_{t}^{2}-\Delta+q_{1}(t,x)}u_{1,\lambda}(t,x)=0  & \mbox{in}\,\,Q, \\

    u_{1,\lambda}(0,x)=u_{2,\lambda}(0,x),\,\,\,\p_{t}u_{1,\lambda}(0,x)=\p_{t}u_{2,\lambda}(0,x) & \mbox{in}\,\,\Omega,\\

  u_{1,\lambda}(t,x)=u_{2,\lambda}(t,x):=f_{\lambda}(t,x), & \mbox{on}\,\,\Sigma.
  \end{array}
\right.
$$
Putting $u_{\lambda}(t,x)=u_{1,\lambda}(t,x)-u_{2,\lambda}(t,x)$, we  get that
$$
\left\{
  \begin{array}{ll}
    \para{\p_{t}^{2}-\Delta+q_{1}(t,x)}u_{\lambda}(t,x)=q(t,x)u_{2,\lambda}(t,x) & \mbox{in}\,\,Q \\
    u_{\lambda}(0,x)=\p_{t}u_{\lambda}(0,x)=0 & \mbox{in}\,\,\Omega \\
    u_{\lambda}(t,x)=0 & \mbox{on}\,\,\Sigma.
  \end{array}
\right.
$$
Applying Lemma \ref{Lem2.2}, once more for $\lambda$  large enough,  we may find a
geometrical optic solution $v_{\lambda}$ to the backward wave equation
$$\para{\p_{t}^{2}-\Delta+q_{1}(t,x)}v_{\lambda}(t,x)=0,\,\,\,\mbox{in}\,\, Q,$$
of the form
\begin{equation}\label{Eq5.38}
v_{\lambda}(t,x)=a(t,x)e^{-i\lambda(x.\omega+t)}+R_{1,\lambda}(t,x),
\end{equation}
where  $R_{1,\lambda}$ satisfies
$$\p_{t}R_{1,\lambda|t=T}=R_{1,\lambda|t=T}=0,\,\,\,\,R_{1,\lambda|\Sigma}=0,$$
and
\begin{equation}\label{Eq5.39}
\|R_{1,\lambda}\|_{L^{2}(Q)}\leq\frac{C}{\lambda}\,\|\varphi\|_{H^{3}(\R^{n})}.
\end{equation}
By integrating by parts and using the Green's formula, one gets
\begin{eqnarray}\label{Eq5.40}
\displaystyle\int_{Q}q(t,x)\!\!\!\!\!\!&\,&\!\!\!\!\!u_{2,\lambda}(t,x)v_{\lambda}(t,x)\,dx\,dt=
\displaystyle\int_{\Sigma} (\mathcal{I}^{1}_{q_{2}}-\mathcal{I}^{1}_{q_{1}})(\psi_{\lambda})v_{\lambda}(t,x)\,d\sigma\,dt\cr
&&+\displaystyle\int_{\Omega}\displaystyle\para{\mathcal{I}_{q_{2}}^{2}-\mathcal{I}_{q_{1}}^{2}}(\psi_{\lambda})\,\p_{t}v_{\lambda}(T,.)\,dx
-\displaystyle\int_{\Omega}\para{\mathcal{I}_{q_{2}}^{3}-\mathcal{I}_{q_{1}}^{3}}(\psi_{\lambda})\,v_{\lambda}(T,.)\,dx,
\end{eqnarray}
where $$\psi_{\lambda}= (u_{2,\lambda|\Sigma},u_{2,\lambda|t=0},\p_{t}u_{2,\lambda|t=0}).$$
Next, we proceed by a similar way as in the proof of Lemma \ref{Lem4.1}, we get
$$\begin{array}{lll}
|\displaystyle\int_{Q}q(t,x)a^{2}(t,x)\,dx\,dt|
&\leq&\|\mathcal{I}_{q_{2}}-\mathcal{I}_{q_{1}}\|\,\|\psi_{\lambda}\|_{H^{1}(\Sigma)\times H^{1}(\Omega)\times L^{2}(\Omega)}\|\phi_{\lambda}\|_{L^{2}(\Sigma)\times L^{2}(\Omega)\times L^{2}(\Omega)}\\
&&+\displaystyle\frac{C}{\lambda}\,\|\varphi\|^{2}_{H^{3}(\R^{n})}
\end{array}$$
where
 $$\phi_{\lambda}=(v_{\lambda|_{\Sigma}},\,v_{\lambda|t=T},\,\p_{t}v_{\lambda|t=T}).$$
Further, as  $R_{i,\lambda|\Sigma}=0$, $\mbox{for}\,\,i=1,2,$, we
deduce that 

$$\begin{array}{lll} 
|\displaystyle\int_{Q}q(t,x)\,a^{2}(t,x)\,dx\,dt|&\leq&
\|\mathcal{I}_{q_{2}}-\mathcal{I}_{q_{1}}\|\,\|\psi_{1,\lambda}\|_{H^{2}(Q)\times H^{1}(\Omega)\times L^{2}(\Omega)}\|\phi_{1,\lambda}\|_{H^{1}(Q)\times L^{2}(\Omega)\times L^{2}(\Omega)}\\
&&+\displaystyle\frac{C}{\lambda}\,\|\varphi\|^{2}_{H^{3}(\R^{n})},
\end{array}$$
where $$\phi_{1,\lambda}=\displaystyle\left(v_{\lambda}-R_{1,\lambda}, \,v_{\lambda|t=T},\,\p_{t}v_{\lambda|t=T}\right),\,\,\,\,
\psi_{1,\lambda}=\para{u_{2,\lambda}-R_{2,\lambda},\,u_{2,\lambda|t=0},\,\p_{t}u_{2,\lambda|t=0}}. $$
Using the fact that $R_{1,\lambda}(T,.)=\p_{t}R_{1,\lambda}(T,.)=0$\, on $\Omega$, we have
   $$\|\phi_{1,\lambda}\|_{H^{1}(Q)\times L^{2}(\Omega)\times L^{2}(\Omega)}\leq
C\lambda\,\|\varphi\|_{H^{3}(\R^{n})},$$
On the other hand, since $R_{2,\lambda}(0,.)=\p_{t}R_{2,\lambda}(0,.)=0$ on $\Omega$, we have
$$\begin{array}{lll}
\|\psi_{1,\lambda}\|_{H^{2}(Q)\times H^{1}(\Omega)\times L^{2}(\Omega)}&\leq& \|u_{2,\lambda}-R_{2,\lambda}\|_{H^{2}(Q)}+\|u_{2,\lambda|t=0}\|_{H^{1}(\Omega)}+\|\p_{t}u_{2,\lambda|t=0}\|_{L^{2}(\Omega)}\\
&\leq&
C \lambda^{2}\,\|\varphi\|_{H^{3}(\R^{n})},
\end{array}$$
Therefore by extending  $q(t,x)$ by zero outside $Q$ and recalling (\ref{Eq2.3}),
we find out that
$$
|\int_{0}^{T}\int_{\R^{n}}q(t,x-t\omega)\,\varphi^{2}(x)\,dx\,dt|\leq C\para{\lambda^{3}\|\mathcal{I}_{q_{2}}-\mathcal{I}_{q_{1}}\|+
\frac{1}{\lambda}}\|\varphi\|^{2}_{H^{3}(\R^{n})}.
$$
Now, in order to complete the proof of Lemma \ref{Lem5.1}, it will be enough to fix $y\in\R^{n}$, consider $(\varphi_{\varepsilon})_{\varepsilon}$ defined as before, and proceed as in the proof of Lemma \ref{Lem3.2}.
By repeating the arguments used in the previous sections, we complete the proof of Theorem \ref{Thm3}.
\end{Demo}

\textbf{Acknowledgements} \\
I would like to thank the professor Mourad Bellassoued for his assistance, for many helpful suggestions he made
and for his careful reading of the manuscript.


\end{document}